\documentclass[reqno,12pt]{amsart}
\usepackage{amscd,amssymb}

\theoremstyle{plain}
\newtheorem{Thm}[subsection]{Theorem}
\newtheorem{Cor}[subsection]{Corollary}
\newtheorem{Lem}[subsection]{Lemma}
\newtheorem{Prop}[subsection]{Proposition}
\newtheorem{Conj}[subsection]{Conjecture}

\theoremstyle{definition}
\newtheorem{Def}[subsection]{Definition}

\theoremstyle{remark}

\newtheorem{Rem}[subsection]{Remark}

\numberwithin{equation}{section}
\renewcommand{\rm}{\normalshape}

\newif\ifShowLabels
\ShowLabelstrue
\newdimen\theight
\def\TeXref#1{%
	\leavevmode\vadjust{\setbox0=\hbox{{\tt
		\quad\quad  {\small \rm #1}}}%
	\theight=\ht0
	\advance\theight by \lineskip
	\kern -\theight \vbox to
	\theight{\rightline{\rlap{\box0}}%
	\vss}%
	}}%

\ShowLabelsfalse

\renewcommand{\sec}[2]{\section{#2}\label{S:#1}%
	\ifShowLabels \TeXref{{S:#1}} \fi}
\newcommand{\ssec}[2]{\subsection{#2}\label{SS:#1}%
	\ifShowLabels \TeXref{{SS:#1}} \fi}

\newcommand{\refss}[1]{Section ~\ref{SS:#1}}
\newcommand{\reft}[1]{Theorem ~\ref{T:#1}}
\newcommand{\refl}[1]{Lemma ~\ref{L:#1}}

\newcommand{\refe}[1]{\eqref{E:#1}}

\newenvironment{thm}[1]%
	{ \begin{Thm} \label{T:#1}  \ifShowLabels \TeXref{T:#1} \fi }%
	{ \end{Thm} }

\renewcommand{\th}[1]{\begin{thm}{#1} \sl }
\renewcommand{\eth}{\end{thm} }

\newenvironment{lemma}[1]%
	{ \begin{Lem} \label{L:#1}  \ifShowLabels \TeXref{L:#1} \fi }%
	{ \end{Lem} }
\newcommand{\lem}[1]{\begin{lemma}{#1} \sl}
\newcommand{\elem}{\end{lemma}}

\newenvironment{propos}[1]%
	{ \begin{Prop} \label{P:#1}  \ifShowLabels \TeXref{P:#1} \fi }%
	{ \end{Prop} }
\newcommand{\prop}[1]{\begin{propos}{#1}\sl }
\newcommand{\eprop}{\end{propos}}

\newenvironment{corol}[1]%
	{ \begin{Cor} \label{C:#1}  \ifShowLabels \TeXref{C:#1} \fi }%
	{ \end{Cor} }
\newcommand{\cor}[1]{\begin{corol}{#1} \sl }
\newcommand{\ecor}{\end{corol}}

\newenvironment{defeni}[1]%
	{ \begin{Def} \label{D:#1}  \ifShowLabels \TeXref{D:#1} \fi }%
	{ \end{Def} }
\newcommand{\defe}[1]{\begin{defeni}{#1} \sl }
\newcommand{\edefe}{\end{defeni}}

\newenvironment{remark}[1]%
	{ \begin{Rem} \label{R:#1}  \ifShowLabels \TeXref{R:#1} \fi }%
	{ \end{Rem} }
\newcommand{\rem}[1]{\begin{remark}{#1}}
\newcommand{\erem}{\end{remark}}

\newcommand{\eq}[1]%
	{ \ifShowLabels \TeXref{E:#1} \fi 
	   \begin{equation} \label{E:#1} }
\newcommand{\eeq}{ \end{equation} }

\newcommand{\prf}{ \begin{proof} }
\newcommand{\eprf}{ \end{proof} }

\newcommand\alp{\alpha}		
\newcommand\bet{\beta}

\newcommand\eps{\varepsilon}		

\newcommand\tet{\theta}		
\newcommand\iot{\iota}

\newcommand\ome{\omega}		\newcommand\Ome{\Omega}

\newcommand\calD{{\mathcal{D}}}

\newcommand\calL{{\mathcal{L}}}

\newcommand\calO{{\mathcal{O}}}

\newcommand\calR{{\mathcal{R}}}
\newcommand\calS{{\mathcal{S}}}

\newcommand\RR{\mathbb{R}}

\newcommand\CC{\mathbb{C}}

\newcommand\nek{,\ldots,}
\newcommand\sdp{\times \hskip -0.3em {\raise 0.3ex
\hbox{$\scriptscriptstyle |$}}} 

\newcommand\Cone{\operatorname{Cone}}

\newcommand\oalp{{\overline{\alpha}}}

\newcommand\tilg{{\widetilde{g}}}

\newcommand\tilH{{\widetilde{H}}}

\newcommand\tilp{{\widetilde{p}}}

\newcommand\tilv{{\widetilde{v}}}

\renewcommand{\>}{\rangle}
\newcommand{\<}{\langle}

\def\o{\calO}
\def\w{\Ome}
\def\sw{\calS'\w}

\def\r{\RR^n}
\def\d{\partial}
\def\b{\bullet}
\def\n{\nabla}  

\def\Cone{\operatorname{Cone}}

\begin{document}


\title[The cohomology of the remote fiber of a real polynomial
map]{Tempered currents and the cohomology of the remote fiber of a
real polynomial map}
\author{Alexander Braverman  \and Maxim Braverman}
\address{School of Mathematical Sciences\\
 	Tel-Aviv University\\
	Ramat-Aviv 69978, Israel}
\email{braval@math.tau.ac.il}
\address{Department of Mathematics\\
        Ohio State University\\
        Columbus, Ohio 43210}
\email{maxim@math.ohio-state.edu}

\begin{abstract}
Let $p:\RR^n\to\RR$ be a polynomial map. Consider the complex
$\sw^{\b}(\RR^n)$ of tempered currents on $\RR^n$ with the twisted
differential $d_p=d-dp$ where $d$ is the usual exterior differential
and $dp$ stands for the exterior multiplication by $dp$. Let $t\in\RR$
and let $F_t=p^{-1}(t)$. 
In this paper we prove that the reduced
cohomology $\tilH^k(F_t;\CC)$ of $F_t$ is isomorphic to
$H^{k+1}(\sw^{\b}(\RR^n),d_p)$ in the case when $p$ is homogeneous and
$t$ is any positive real number. 
We conjecture that this isomorphism holds for any polynomial $p$, 
for $t$ large enough (we call the $F_t$ for $t\gg 0$ {\it the remote fiber}
of $p$) and we prove this conjecture for polynomials that satisfy
certain technical condition (cf. \reft{main'}). 
The result is analogous to that of
A.~Dimca and M.~Saito (\cite{DimSa}), who give a similar (algebraic) way
to compute the reduced cohomology of the generic fiber of a complex
polynomial.
\end{abstract}
\maketitle

\sec{introd}{Introduction}

\ssec{0.1}{The Dimca-Saito theorem}Let $p:\CC^n \to \CC$ be a complex
polynomial. Let $F$ denote the {\em generic fiber} of $p$ (it is well-defined
as a topological space). In \cite{DimSa}, A.~Dimca and M.~Saito have given the
following algebraic way to compute the cohomology of $F$. Let
$\Ome^{\b}$ denote the De Rham algebra of polynomial differential
forms on $\CC^n$. Define a differential $d_p$ on $\Ome^{\b}$ by
$$
	d_p(\ome)=d\ome-dp\wedge \ome
$$
\begin{Thm}[\textbf{Dimca-Saito}]\label{T:DS}\sl
  There exists an isomorphism 
  $$
    H^{k+1}(\Ome^{\b},d_p)\simeq \tilH^{k}(F,\CC)
			\quad\text{for }k=0,1\nek n-1
  $$
  where $\tilH^{\b}(F,\CC)$ denotes the reduced cohomology of $F$ with 
  coefficients in $\CC$.
\end{Thm}
\ssec{0.2}{The main result} The main purpose of this paper is to describe
certain real analogue of \reft{DS}. Namely, let now $p:\RR^n\to \RR$ be
a real polynomial. Then (cf. \cite{Varch1}) the topological type of
the fiber $F_t=p^{-1}(t)$ does not depend upon $t$ provided $t$ is
large enough. We shall refer to $F_t$ as to {\it remote fiber} of $p$
and we shall be interested in the cohomology of $F_t$.

Let $\o:=\CC\, [x_1\nek x_n]$ denote the ring of complex polynomials on
$\RR^n$ and let \/ $\calS'(\RR^n)$ \/ denote the $\o$ module of tempered
complex valued distributions on $\RR^n$.  Let $\w^\b$ be the complex
of global algebraic differential forms on $\RR^n$. Consider the space
$$
        \sw^\b(\r)=\calS'(\RR^n)\otimes_{\o} \Ome^\b
$$ 
of tempered currents  and define a differential $d_p:\sw^{\b}(\r)\to
\sw^{\b+1}(\r)$ on $\sw^{\b}(\r)$ by
$$
        d_p(\ome)=d\ome - dp\wedge \ome \qquad 
                                \text{for} \quad \ome\in\sw^{\b}(\r).
$$

In this paper we discuss the following
\begin{Conj}\label{T:mainC}\sl
   For any real polynomial  $p:\RR^n\to \RR$ the following
   isomorphism holds:
   \eq{main}
        H^{k+1}(\sw^{\b}(\RR^n),d_p)=\tilH^{k}(F_t;\CC), \qquad k=0,1,\dots,n-1
   \end{equation}
   where $\tilH$ denotes reduced cohomology. 
\end{Conj}

In particular, we prove the following


\th{main} Assume that $p:\r\to \RR$ is a homogeneous polynomial map of degree
   $m$, i.e.  $p(sx)=s^mp(x)$ for any $x\in\r, \ s\in \RR$. For any
   $t>0$, the isomorphism \refe{main} holds.
\eth

\rem{rem}
{\it a)} \ By definition, the reduced cohomology of any topological
space $X$ is the cohomology of the complex
$$\begin{CD}
  0\to \CC @>\eps>> H^0(X;\CC)@>0>>\cdots @>0>>H^k(X;\CC)\to\cdots
\end{CD}$$ 
where $\eps=0$ if $X$ is empty and $\eps$ is the natural map
$\CC=H^0(pt;\CC)\to H^0(X;\CC)$ coming from the projection $X\to pt$
otherwise. Therefore in the case when $X$ is empty
one should have $\tilH^{-1}(X;\CC)=\CC$ (but $\tilH^{-1}(X;\CC)=0$ if
$X$ contains at least one point). With this convention, \reft{main}
remains true also for $k=-1$.

{\it b) \ }Note that if instead of the complex $\sw^{\b}(\r)$ we
considered the complex of all currents with the same differential
$d_p$, then we would get a complex quasi-isomorphic to the usual
complex of currents on $\RR^n$ with the ordinary exterior differential
$d$ (since $d$ and $d_p$ are conjugate to one another by means of the
function $e^p$).  Therefore if we do not impose any growth conditions
on our currents we will not get any interesting cohomology.
\erem

\ssec{0.4}{Sketch of the proof} \ 
{\it Step 1. \ }Let $U_t=\{x\in\RR^n:\ p(x)>t\}$ (note that $U_t$ might
be empty).  Then,  $U_t$ is diffeomorphic to the product $F_t\times
(0,\infty)$, for any $t>0$.   Using the long cohomological
sequence of the pair $(\RR^n,U_t)$ one can easily see that
$\tilH^{k-1}(F_t,\CC)=H^k(\RR^n,U_t;\CC)$.

{\it Step 2. \ }Let $\calD'\Ome^{\b}(U_t)$ denote the complex of all
currents on $U_t$. In \refss{S'} we define certain subcomplex
$\sw^{\b}(U_t)$ of $\calD'\Ome^{\b}(U_t)$ ({\it the complex of
tempered currents on $U_t$}) and prove that its natural inclusion into
$\calD'\Ome^{\b}(U_t)$ is a quasi-isomorphism.

{\it Step 3. \ }Let $\calD'\Ome^{\b}(\RR^n)$ denote the complex of all
currents on $\RR^n$.  Let $\theta (s)$ be a smooth function on $\RR$,
such that $\theta (s)=s$ for $s<1$ and $\theta (s)=0$ for
$s>2$. Define $\tilp:\RR^n\to \RR$ by $\tilp(x)=\theta (p(x))$. Let
$\calS'_p\Ome^{\b}(\RR^n)$ denote the space of all currents $\ome$ on
$\RR^n$, such that $e^{\tilp}\ome\in \sw^{\b}(\RR^n)$.  Then we show in
\refl{S'p} that $\calS'_p\Ome^{\b}(\RR^n)$ is a subcomplex of
$\calD'\Ome^{\b}(\RR^n)$ and the natural embedding
$\calS'_p\Ome^{\b}(\RR^n)\hookrightarrow\calD'\Ome^{\b}(\RR^n)$ is a
quasi-isomorphism.

Let now $\rho$ denote the natural map from $\calS'_p\Ome^{\b}(\RR^n)$ to
$\sw^{\b}(U_t)$ (restriction to $U_t$). It follows from step 2 and from 
the above statement 
that the complex 
$$
	\Cone^\b(\rho)= \calS'_p\Ome^\b(\r)\oplus
		\calS'\Ome^{\b-1}(U_t)
$$
computes the relative cohomology $H^{\b}(\RR^n,U_t;\CC)$.

{\it Step 4. \ }The map $\Phi_1:\ \ome\to e^{-p}\ome$ defines
a morphism of complexes $\sw^{\b}(\RR^n)\to \calS'_p\Ome^\b(\r)$. 
Moreover, every element
in the image of $\Phi_1$ is rapidly decreasing along the rays 
$R_x=\big\{ sx:\, s>0\big\}$, for any $x\in F_t$. 
 This  enables us 
to extend $\Phi_1$ to an explicit map 
$
\Phi:\ \sw^{\b}(\RR^n)\to \Cone^{\b}(\rho)
$.  
In order to do that we need the following notations.

Let $\mu_s:\r\to \r$ denote the multiplication by $s$ and let
$\mu_s^*:\calD'\Ome^\b(\r)\to \calD'\Ome^\b(\r)$ be the corresponding
pull-back map.

Consider the {\em Euler vector field} 
$
	\calR=\sum_{i=1}^n  x_i\frac \d{\d x_i}
$
on $\r$ and let $\iot_\calR, \ \calL_\calR$ denote the interior
multiplication by $\calR$ and the Lie derivative along $\calR$. Then
\eq{dmu'}
	\frac d{ds}\mu_s^*(\ome)= \mu_s^*\big(\calL_\calR\ome\big)s^{-1}
		\quad \text{for any} \quad\ome\in \calD'\Ome^\b(\r).
\end{equation}

We define the map 
$
	\Phi: \calS'\Ome^\b(\r)\to \Cone^\b(\rho)=
		\calS'_p\Ome^\b(\r)\oplus \calS'\Ome^{\b-1}(U_t)
$ 
by the formula
\eq{2.4'}
	\Phi:\ome\mapsto \big(\Phi_1\ome, \Phi_2\ome\big)=
	  \left(e^{-p}\ome,\, -\int_1^\infty
	  \mu_s^*(e^{-p}\iot_{\calR}\ome)\, \frac{ds}s\right).
\end{equation}
The integral in \refe{2.4'} converges since $e^{-p(sx)}$ decreases
exponentially in $s$ as $s$ tends to infinity.  One uses \refe{2.4'}
to show that $\Phi$ commutes with differentials.

Finally we prove by an explicit calculation that $\Phi$ is a
quasi-isomorphism.  Therefore $\sw^{\b}(\RR^n)$ computes
$H^{\b}(\RR^n,U_t;\CC)$, which is isomorphic to
$\tilH^{\bullet-1}(F_t,\CC)$ by step 1.

\ssec{general}{The general case} Let now $p:\RR^n\to \RR$ be an
arbitrary polynomial.  Set $v=\frac{\n p}{|\n p|^2}$. Then the Lie
derivative of $p$ along $v$ is equal to 1. In the Appendix we show that
the flow of $v$ is globally defined on $U_t$ if $t$ is large
enough. We denote this flow by $g_s:U_t\to U_t$ and let $g_s^*:
\calD'\Ome^\b(U_t)\to \calD'\Ome^\b(U_t)$ be the corresponding
pull-back of currents. Then $g_s^*(p)=p+s$. In particular, we obtain a
new proof of topological equivalence of the fibers $F_t$ with $t\gg
0$.

Denote 
$$
	\tilv \ = \ pv
$$
and let $\iot_\tilv, \ \calL_\tilv$ denote the interior multiplication
by $\tilv$ and the Lie derivative along $\tilv$. The flow $\tilg_s$ of
$\tilv$ is defined on $U_t, t\gg 0$. The flow $\tilg_s$ and the vector
field $\tilv$ are connected by the formula
$$
	\frac d{ds}\tilg_s^*(\ome) \ = \ 
		\tilg_s^*\big(\calL_\tilv\ome\big)
		\quad \text{for any} \quad\ome\in \calD'\Ome^\b(\r),
$$
which is similar to \refe{dmu'} (if $p$ is a homogeneous polynomial of
degree $m$ then $\mu_{s}= g_{m\ln s}$). One can easily check that
$\tilg_s^*(p)=e^sp$.

One can try to define a map  $\Phi: \calS'\Ome^\b(\r)\to
\Cone^\b(\rho)$ by formula 
$$
	\Phi:\ome\mapsto \big(\Phi_1\ome, \Phi_2\ome\big)=
	  \left(e^{-p}\ome,\, -\int_1^\infty
	  \tilg_s^*(e^{-p}\iot_{\tilv}\ome)\, ds\right),
$$
similar to \refe{2.4'}. The only problem here is that
we were not able to prove that the integral in the definition of
$\Phi_2$ converges to a tempered current. However, if the map 
$\Phi_2:\calS'\Ome^\b(\r)\to \calS'\Ome^{\b-1}(\r)$ is well defined a
verbatim repetition of the proof of \reft{main} gives the following
\th{main'}Suppose that $p:\RR^n\to \RR$ is a real polynomial and
   $\tilg_s, s>0$ is a one-parameter semigroup of diffeomorphisms
   $U_t\to U_t$ such that $\tilg_s^*(p)=e^{ms}p$. Let
   $\tilv=\frac{d}{ds}_{|_{s=0}}g_s$.  If for any tempered current
   $\ome$ the integral
   $$
	\int_1^\infty
	  \tilg_s^*(e^{-p}\iot_{\tilv}\ome)\, ds
   $$ 
   converges to a tempered current, then the isomorphism \refe{main}
   holds.
\eth

\ssec{example}{Example} Consider the polynomial of two variables
$p(x,y)=x^2-x-y$. Set $U_0=\big\{ (x,y)\in \RR^2:\, p(x,y)>0\big\}$
and define a one parameter semigroup $g_s$ of diffeomorphisms of $U_0$
by the formula
$$
	g_s(x,y) \ = \ \left( e^{s/2}x, e^{s}x-e^{s/2}x+e^{s}y \right).
$$
Then $g_s^*p=e^sp$. Clearly, all other conditions of \reft{main'} are
satisfied. Hence, the isomorphism \refe{main} holds for $p(x,y)$.

\subsection*{Acknowledgments}It is a great pleasure for us to express
our gratitude to J.~Bernstein and M.~Farber; the paper was
considerably influenced by communications with them. It was M.~Farber
who sugested to use the map \refe{2.4'} for the study of $H^\b(\Ome^\b,d_p)$.

We are also
thankful to N.~Zobin, M.~Zaidenberg and S.~Kaliman. 

The first author would like to thank Institute for Advance Study for
hospitality.
\sec{cur}{Complexes of Currents}

In this section we review some facts about complexes of currents which
will be used in the proof of \reft{main}.

Let $p:\r\to \RR$ be a homogeneous polynomial map of degree $m$,
i.e. $p(sx)=s^mp(x)$. Let $U_t=\nolinebreak{\big\{x\in\r:\,
p(x)>t\big\}}$, where $t\in \RR$.

\ssec{D'}{The complex of currents} By
$\Ome_c^\b(U_t)$ we denote the De Rham complex of compactly supported
complex valued $C^\infty$-forms on $U_t$. The cohomology of
$\Ome_c^\b(U_t)$ is called the {\em compactly supported cohomology} of
$U_t$.

Recall that if 
$
	0\to C^0\overset{d}{\to} C^1\overset{d}{\to}
			\cdots\overset{d}{\to}  C^n\to 0 
$ 
is a complex of topological vector spaces then the dual complex to
$(C^\b,d)$ is, by definition, the complex
$$\begin{CD}
	0\to (C^n)^*@>d^*>> (C^{n-1})^*@>d^*>>\cdots @>d^*>> (C^0)^*\to 0,
\end{CD}$$
where $(C^i)^*$ denotes the topological dual of the space $C^i$ and
$d^*$ denotes the  adjoint operator of $d$.

The complex of currents $\calD'\Ome^\b(U_t)$ on $U_t$ is the complex
dual to $\Ome_c^\b(U_t)$.  By the Poincar\'e duality for non-compact
manifolds (cf. \cite{BottTu}), the cohomology of $\calD'\Ome^\b(U_t)$
is equal to the cohomology of $U_t$.

Analogously, one defines the complex $\calD'\Ome^\b(\r)$ of currents
on $\r$.  

Let $r:\calD'\Ome^\b(\r)\to \calD'\Ome^\b(U_t)$ be the
restriction. Recall that the cone $\Cone^\b(r)$ of $r$ is the complex
$$
	\Cone^\b(r)=\calD'\Ome^\b(\r)\oplus \calD'\Ome^{\b-1}(U_t), \qquad
			d:(\ome,\alp)\mapsto (d\ome,\ome-d\alp). 
$$
The cohomology of $\Cone^\b(r)$ is equal to the relative cohomology
$H^\b(\r,V;\CC)$ of the pair $(\r,V)$.

\ssec{S'}{The complex of tempered currents} The space
$\calS(\r)$ of  Schwartz (rapidly decreasing) functions on $\r$ is the set
of all $\phi\in C^\infty(\r)$ such that for any linear differential operator
$L:C^\infty(\r)\to  C^\infty(\r)$ with polynomial coefficients
\eq{1.1}
	\sup_{x\in\r} |L\phi(x)| <\infty.
\end{equation}
The topology in $\calS(\r)$ defined by the semi-norms in the left-hand
side of \refe{1.1} makes $S(\r)$ a Fr\'echet space.

Recall that by $\Ome^\b$ we denote the De Rham complex of
global algebraic differential forms on $\r$. The  {\em complex of Schwartz
forms} on $\r$ is the complex
$$
	\calS\Ome^\b(\r)=\calS(\r)\otimes_{\o} \Ome^\b  \qquad
$$ 
with natural differential. By the {\em complex of Schwartz forms on $U_t$}
we will understand the subcomplex of $\calS\Ome^\b(\r)$ consisting of the
forms $\ome$ such that there exists a real number $\eps=\eps(\ome)$
such that the support of $\ome$ lies in $U_{t+\eps}$.

The complex $\calS'\Ome^\b(\r)$ of {\em tempered currents} on $\r$ is, by
definition, the dual complex to $\calS\Ome^\b(\r)$. Similarly,
the complex $\calS'\Ome^\b(U_t)$ of {\it tempered currents on $U_t$}
is, the dual  complex to $\calS\Ome^\b(U_t)$. 

\lem{t1t2} For any $t_1> t_2 >0$ the natural map
  $i: \calS'\Ome^\b(U_{t_2})\to
  \calS'\Ome^\b(U_{t_1})$ is a homotopy equivalence of complexes.
\elem
\prf
Let $\mu_s:\r\to \r$ denote the multiplication by $s$ and let
$\mu_s^*:\calD'\Ome^\b(\r)\to \calD'\Ome^\b(\r)$ be the corresponding
pull-back map. Clearly, $\mu_s^*$ preserves the space of tempered currents.

Set $\tau=(t_1/t_2)^{1/m}$. Then $\mu_\tau(U_{t_2})=U_{t_1}$.  In
particular, we can consider $\mu^*_{\tau}$ as a map from
$\calS\Ome^\b(U_{t_1})$ to $\calS\Ome^\b(U_{t_2})$. To prove the lemma
we will show that $\mu^*_\tau$ is  a homotopy inverse of $i$.

Consider the {\em Euler vector field} 
$$
	\calR=\sum_{i=1}^n  x_i\frac \d{\d x_i}
$$
on $\r$ and let $\iot_\calR, \ \calL_\calR$ denote the interior
multiplication by $\calR$ and the Lie derivative along $\calR$. Then
\eq{dmu}
	\frac d{ds}\mu_s^*(\ome)= \mu_s^*\big(\calL_\calR\ome\big)s^{-1}
		\quad \text{for any} \quad\ome\in \calD'\Ome^\b(\r).
\end{equation}
Note that if $\ome$ is a tempered current so are $\iot_\calR\ome$ and
$\calL_\calR\ome$. 

For any current $\ome$, set
$$
	H\ome \ = \ 
	    \int^{\tau}_1 \mu^*_s(\iot_{\calR}\ome)\, \frac{ds}s.
$$
The operators $\mu_s^*$ and $\iot_{\calR}$
preserve the space of tempered currents. Hence, so does $H$. Using
\refe{dmu} and the Cartan homotopy formula 
\eq{Cartan}
	\calL_\calR=d\iot_\calR+ \iot_\calR d
\end{equation}
we obtain 
$$
	(dH+Hd)\, \ome \ = \ \mu_\tau^*\ome \ - \ \ome,
$$
for any current $\ome$. The lemma is proven.
\eprf


\lem{S'D'} For any $t> 0$ the embedding 
  $\calS'\Ome^\b(U_t)\hookrightarrow \calD'\Ome^\b(U_t)$ is a homotopy
  equivalence of complexes.  In particular, the cohomology of the
  complex $\calS'\Ome^\b(U_t)$ is equal to the cohomology
  $H^\b(U_t;\CC)$ of $U_t$.
\elem 
\prf 
By \refl{t1t2} it is enough to show that the embedding 
$i:\calS'\Ome^\b(U_1)\hookrightarrow \calD'\Ome^\b(U_1)$ is a
quasi-isomorphism.

Let $h_s:U_1\to U_1, \ s>0$ denote the map defined by the formula
$$
	h_s:x\mapsto \frac{1+s}{1+s|x|}\cdot x.
$$
Here $|x|$ denotes the norm of the vector $x\in \r$. 
Let $h_s^*:\calD'\Ome^\b(U_1)\to \calD'\Ome^\b(U_1)$ denote the
corresponding pull-back. Then (cf. \refe{dmu})
\eq{dh}
	\frac{d}{ds}h^*_s(\ome) \ = \ 
		  \frac{1-|x|}{(1+s)(1+s|x|)} \, h^*_s(\calL_\calR\ome),
\end{equation}
for any current $\ome$.  Note also that $h_s^*$ preserves the space of
tempered currents.

Clearly,  $h_0$ is the identity map.
The image of $h_1:U_1\to U_1$ lies in the compact set 
$$
	\big\{x\in \RR^n:\ |x|\le 2\big\}.
$$
Hence, $h_1^*\ome$ is a tempered current for any
$\ome\in\calD'\Ome^\b(\r)$. Note also that $h_s^*$ preserves the space of
tempered currents for any $s>0$. We will prove that the map $h_1^*:
\calD'\Ome^\b(U_1)\to \calS'\Ome^\b(U_1)$ is a homotopy inverse of the
embedding $i:\calS'\Ome^\b(U_1)\hookrightarrow \calD'\Ome^\b(U_1)$.

For any $\ome\in \calD'\Ome^\b(U_1)$, set
$$
	H\ome \ = \ 
	   \int_0^1\, 
	\frac{1-|x|}{(1+s)(1+s|x|)} \, h_s^*(\iot_\calR\ome)\, ds. 
$$
Here, $\iot_\calR$ denote the operator of interior multiplication by
$\calR$. 

Using \refe{dh} and \refe{Cartan} we obtain
\eq{homot'}
	(dH+Hd)\, \ome = h_1^*\, \ome -\ome, \qquad 
			\ome\in  \calD'\Ome^\b(U_1).
\end{equation}
Thus the map $i\circ h_1^*: \calD'\Ome^\b(U_1)\to \calD'\Ome^\b(U_1)$
homotopic to the identity map.

Since the operators $h_s^*$ and $\iot_\calR$ preserve the space of
tempered currents, so does $H$. Hence, \refe{homot'} implies that the
map $h_1^*\circ i: \calS'\Ome^\b(U_1)\to \calS'\Ome^\b(U_1)$ is also
homotopic to the identity map.  
\eprf

\ssec{1.5.}{The complex $\calS'_p\Ome^\b(\RR^n)$}
We will need the following twisted version of the complex of tempered
currents on $\r$.

Fix a smooth function $\tet:\RR\to \RR$ such that
$$
	\tet(s)=\cases
			s \quad \text{if} \quad s<1,\\
			0 \quad \text{if} \quad s>2.
		\endcases
$$
and define $\tilp(x)=\tet(p(x)), x\in \RR^n$.  Note that the current
$d\tilp\wedge \ome$ is tempered for any tempered current $\ome$.
\lem{S'p}The space
  $\begin{CD}
	\calS'_p\Ome^\b(\r)=\big\{ \ome \in \calD'\Ome^\b(\r):\, 
				e^{\tilp}\ome\in \calS'\Ome^\b(\r) \big\}
  \end{CD}$
  is a subcomplex of $\calD'\Ome^\b(\r)$ and the embedding
  $\calS'_p\Ome^\b(\r)\hookrightarrow \calD'\Ome^\b(\r)$ is a
  quasi-isomorphism.  
\elem
\prf
Suppose $\ome\in \calS'_p\Ome^\b(\r)$, i.e. $e^{\tilp}\ome\in
\calS'\Ome^\b(\r)$. Then
$$
	e^{\tilp} d\ome=d(e^{\tilp}\ome)-d\tilp\wedge e^{-\tilp}\ome 
							\in \calS'\Ome^\b(\r),
$$
i.e. $d\ome\in \calS'_p\Ome^\b(\r)$. Hence $\calS'_p\Ome^\b(\r)$ is a
subcomplex of $\calD'\Ome^\b(\r)$.

Clearly, the embedding $\calS'_p\Ome^\b(\r)\hookrightarrow
\calD'\Ome^\b(\r)$ induces an isomorphism of 0-cohomology. To prove
\refl{S'p} it remains to show that the $k$-th cohomology
$H^k(\calS'_p\Ome^\b(\r)), \ k>0$ of $\calS'_p\Ome^\b(\r)$ vanishes.

We will use the notation introduced in the proof of \refl{t1t2}. In
particular, $\mu_s:\r\to \r$ is the multiplication by $s$ and
$\calR$ is the Euler vector field on $\r$. 

Let $\ome$ be a closed current, $d\ome=0$. Using \refe{dmu} and the
Cartan homotopy formula  $\calL_\calR=d\iot_\calR+\iot_\calR d$, we obtain
$$
	\ome-\mu_0^*(\ome) =d\, \int_0^1\mu_s^*\big(\iot_\calR\ome\big)\,
							\frac{ds}s.
$$
(Note that the integral in the left hand side converges, because
$\calR$ vanishes at 0). If $\ome$ is a $k$-current, $k>0$, then
$\mu_0^*(\ome)=0$. Hence,
to finish the proof we need only to show that 
\eq{temp}
	\int_0^1\mu_s^*\big(\iot_\calR\ome\big)\, \frac{ds}s\in
		\calS'_p\Ome^\b(\r)
\end{equation}
for any $\ome\in \calS'_p\Ome^\b(\r)$.

Set $\bet =e^\tilp\ome\in \calS'\Ome^\b(\r)$. Then
\eq{new}
	e^\tilp\int_0^1\mu_s^*\big(\iot_\calR\ome\big)\,
					\frac{ds}s=
	\int_0^1 e^{\tilp(x)-\tilp(sx)} \mu_s^*\big(\iot_\calR\bet\big)\,
							\frac{ds}s.
\end{equation}
Since, for any $s\in [0,1]$, the function $\tilp(x)-\tilp(sx)$ is
bounded from above, all the derivatives of the function
$e^{\tilp(x)-\tilp(sx)}$ are bounded by polynomials.  It follows that
$s\mapsto e^{\tilp(x)-\tilp(sx)}
\mu_s^*\big(\iot_\calR\bet\big)s^{-1}$ defines a continuous map
$[0,1]\to \calS'\Ome^\b(\r)$. Hence the current \refe{new} is tempered
and \refe{temp} holds.
\eprf

\ssec{1.7}{}
For any $\ome\in \calS'_p\Ome^\b(\r)$, $t>0$ the restriction of
$\ome$ on $U_t$ is a tempered current on $U_t$. Hence, the restriction
map $\rho:\calS'_p\Ome^\b(\r)\to \calS'\Ome^\b(U_t)$ is defined.
%
\lem{Cone} The complexes $\Cone^\b(\rho)$ and $\Cone^\b(r)$
  (cf. \refss{S'}) are quasi-~isomorphic. In particular,  
  the cohomology of  $\Cone^\b(\rho)$ equals the
  relative  cohomology of the pair $(\r,U_t)$.
\elem
\prf
Let $i:\calS'_p\Ome^\b(\r)\to \calD'\Ome^\b(\r),\, j:\calS'\Ome^\b(U_t)\to
\calD'\Ome^\b(U_t)$ denote the natural inclusions.
Consider the commutative diagram
\eq{1.5}
 \begin{CD}
	\calS'_p\Ome^\b(\r)  @>\rho>>   \calS'\Ome^\b(U_t)\\
	@ViVV				@VVjV  \\
	\calD'\Ome^\b(\r) @>r>>     \calD'\Ome^\b(U_t) 
 \end{CD}
\end{equation}
According to Lemmas \ref{L:S'D'} and \ref{L:S'p} the vertical arrows of this
diagram are quasi-isomorphisms. Hence, \refe{1.5}
induces a quasi-isomorphism between $\Cone^\b(\rho)$ and $\Cone^\b(r)$.   
\eprf

In the proof of \refl{2.6} we will also need the following
\lem{Cone'} Suppose $t_1>t_2>0$ and let 
  $$
	\rho_1:\calS'_p\Ome^\b(\r)\to \calS'\Ome^\b(U_{t_1}),
	\qquad
	 \rho_2:\calS'_p\Ome^\b(\r)\to \calS'\Ome^\b(U_{t_2})
  $$
  denote the corresponding restrictions. The natural map
  $\Cone^\b(\rho_2)\to \Cone^\b(\rho_1)$ is a quasi-isomorphism.
\elem
\prf
  Consider the commutative diagram
\eq{1.6}
 \begin{CD}
	\calS'_p\Ome^\b(\r)  @>\rho_2>>   \calS'\Ome^\b(U_{t_2})\\
	@|				@VVV  \\
	\calS'_p\Ome^\b(\r) @>\rho_1>>     \calS'\Ome^\b(U_{t_1}) 
 \end{CD}
\end{equation}
By \refl{t1t2}, the right vertical arrow of this diagram is a
quasi-isomorphism. Hence, \refe{1.6} induces a quasi-isomorphism
between $\Cone^\b(\rho_2)$ and $\Cone^\b(\rho_1)$.
\eprf

\sec{proof}{Proof of Theorem 0.2}

\ssec{FtUt}{Cohomology of $F_t$ as relative cohomology} Fix $t> 0$ and
set  
$$
	F_t=p^{-1}(t); \qquad U_t= \big\{x\in\RR^n:\, p(x)>t\big\}.
$$ 
Then $U_t$ is diffeomorphic to the the product $F_t\times
(0,+\infty)$. In particular, $U_t$ has the same cohomology as
$F_t$. Using the long exact sequence of the pair $(\r,U_t)$, we obtain
\eq{2.3}
	 \tilH^k(F_t;\CC)=H^{k+1}(\r,U_t;\CC),  \qquad k=0,1,\dots, n-1,
\end{equation}
where $H^\b(\r,U_t;\CC)$ denotes the relative cohomology of the pair
$(\r,U_t)$ and $\tilH^\b(F_t;\CC)$ denotes the reduced cohomology of
$F_t$.  

From \refl{Cone} and \refe{2.3} we see that to prove \reft{main} it
is enough to show that the complexes $(\calS'\Ome^\b(\r),d_p)$ and
$\Cone^\b(\rho)$ are quasi-isomorphic. 

\ssec{S'Cone}{A map from $\calS'\Ome^\b(\r)$ to $\Cone^\b(\rho)$} 
Recall that by $\mu_s:\RR^n\to \RR^n$ we denote the multiplication by
$s\in\RR$. Then $\mu_s(U_t)=U_{s^mt}$. In particular, if $s\ge 1$, then
$\mu_s$ may be considered as a map from $U_t$ to itself. 
Let $\mu^*_s:\calD'\Ome^\b(U_t)\to \calD'\Ome^\b(U_t)$ denote the
corresponding pull-back map. Then 
$\mu_s^*p(x)=p(\mu_sx)=s^mp(x)$.

Recall also that $\iot_{\calR}$ denote the interior multiplication by
the Euler vector field $\calR=\sum x_i\frac{\d}{\d x_i}$. Note that if
$\ome$ is a tempered current on $U_t$ than so are $\iot_{\calR}\ome$ and
$\mu_s^*\ome$.

Recall from \refss{1.7} that $\rho: \calS'_p\Ome^\b(\r)\to
\calS'\Ome^\b(U_t)$ denotes the  restriction.
We define the map 
$$
	\Phi: \calS'\Ome^\b(\r)\to \Cone^\b(\rho)=
		\calS'_p\Ome^\b(\r)\oplus \calS'\Ome^{\b-1}(U_t)
$$ 
by the formula
\eq{2.4}
	\Phi:\ome\mapsto \big(\Phi_1\ome, \Phi_2\ome\big)=
	  \left(e^{-p}\ome,\,  
		-\int_1^\infty 
	\mu_s^*(e^{-p}\iot_{\calR}\ome)\, \frac{ds}s\right).
\end{equation}
The integral in \refe{2.4} converges since $e^{-p(sx)}$ decreases
exponentially in $s$ as $s$ tends to infinity. It follows from
\refe{dmu}, \refe{Cartan} that the map $\Phi: \calS'\Ome^\b(\r)\to
\Cone^\b(\rho)$ commutes with differentials, i.e.
$$
	\Phi_1d_p\ome \ = \ d\Phi_1\ome; \qquad \Phi_2
		d_p\ome \ =\ \Phi_1\ome_{|_U} \ - \ d\Phi_2\ome.
$$

\lem{2.5} The map $ H^\b\big( \calS'\Ome^\b(\r),d_p\big)\to
  H^\b\big(\Cone^\b(\rho)\big) $ induced by $\Phi$ is injective.
\elem 
\prf 
Suppose that $\ome$ is a tempered current on $\r$ and that $\Phi\ome$
is a coboundary in $\Cone^\b(\rho)$. Then there exists $\alp\in
\calS'_p\Ome^\b(\r)$ and $\bet\in \calS'\Ome^{\b-1}(U_t)$ such that
$$
	e^{-p}\ome \ = \ d\alp; \qquad 
	   -\int_1^\infty \mu_s^*(e^{-p}\iot_{\calR}\ome)\, \frac{ds}s
			\ = \
				\alp_{|_U} \ - \ d\bet.
$$

Choose $j\in C^\infty(\RR)$ such that $j(s)=0$ if $s\le t+1$ and
$j(s)=1$ if $s\ge t+2$ and set $\phi(x)=j(p(x)), \ (x\in\r)$. Then the
support of $\phi$ is contained in $U_t$. Hence $\phi\bet$ may be considered
as a current on $\r$. Since all the derivatives of $\phi$ are bounded
by polynomials, $\phi\bet\in \calS'_p\Ome^\b(\r)$.

Define $\oalp=\alp-d(\phi\bet)$. Then $\ome= d_p(e^{p}\oalp)$. So to
prove the lemma we only need  to show that $e^p\oalp$ is a tempered
current.  

Let $\psi(x)=j(p(x)-1)$. It is enough to prove that
$e^p(1-\psi)\oalp$ and $e^p\psi\oalp$ are tempered currents.  

Since $(1-\psi)\oalp\in \calS'_p\Ome^\b(\r)$ vanishes
when $p(x)>t+2$, we see from the definition of $\calS'_p\Ome^\b(\r)$ that 
$e^p(1-\psi)\oalp$ is a tempered current. 

On the support of $\psi$ the function $\phi$ is identically equal to
1. Hence $\psi\oalp=\psi(\alp-d\bet)$ and 
\begin{multline}\notag
	e^p\psi\oalp=e^p\psi(\alp-d\bet)= 
	  -\psi\int_1^\infty e^p \mu_s^*(e^{-p}\iot_{\calR}\ome)\, 
							\frac{ds}s=\\
      		-\psi\int_1^\infty e^{p(x)-p(\mu_sx)} 
				\mu_s^*(\iot_{\calR}\ome)\, \frac{ds}s
						\in\calS'\Ome^\b(\r).
\end{multline}

\eprf

\lem{2.6} The map
  $
	H^\b\big( \calS'\Ome^\b(\r),d_p\big)\to H^\b\big(\Cone^\b(\rho)\big)
  $
  induced by $\Phi$ is surjective.
\elem
\prf
Choose $\eps>0$ such that $t-\eps>0$ and set 
$U_{t-\eps}=\big\{x\in\RR^n:\, p(x)>t-\eps\big\}$. 
Let $\rho_\eps:\calS'_p\Ome^\b(\r)\to \calS'\Ome^\b(U_{t-\eps})$ be the
restriction. By \refl{Cone'}, any cohomology class $\xi$ of the
complex $\Cone^\b(\rho)$ may be represented by a pair $(\alp,
\bet_{|_U})$, where $\alp\in \calS'_p\Ome^\b(\r)$ and $\bet\in
\calS'\Ome^\b(U_{t-\eps})$. 

Fix $j\in C^\infty(\RR)$ such that $j(s)=0$ if $s\le t-\eps/2$ and $j(s)=1$
if $s\ge t$ and set $\phi(x)=j(p(x)), \ (x\in\r)$. Then all the
derivatives of $\phi$ are bounded by polynomials and the support of
$\phi$ is contained in $U_{t-\eps}$. Hence, $\phi\bet$, considered as a
current on $\r$, belongs to the space $\calS'_p\Ome^\b(\r)$.

The cohomology class of the pair 
$
	\big(\alp- d(\phi\bet),\, 0\big)\in \Cone^\b(\rho)
$ 
equals $\xi$. Set $\ome=e^p(\alp-d(\phi\bet))$. Since the current 
$\alp-d(\phi\bet)$ vanishes on $U_t$, we see from the definition of the
space $\calS'_p\Ome^\b(\r)$ that $\ome$ is a tempered current.
Clearly, $\Phi\ome=\big(\alp- d(\phi\bet),\,  0\big)$. Hence, $\xi$ belongs to
the image of the map $H^\b\big( \calS'\Ome^\b(\r),d_p\big)\to 
H^\b\big(\Cone^\b(\rho)\big)$
\eprf

From Lemmas \ref{L:2.5} and \ref{L:2.6}, we see that the complexes
$(\calS'\Ome^\b(\r),d_p)$ and $\Cone^\b(\rho)$ are
quasi-isomorphic. \reft{main} follows now from \refl{Cone} and \refe{2.3}.


\appendix
\sec{appendix}{Gradient vector field near infinity}

In this appendix we show that on the set $U_T, T\gg 0$ the gradient
vector field $\n p$ of $p$ is bounded from bellow by $c/|x|$. That
means that the vector field $v=\frac{\n p}{|\n p|^2}$ defined in
\refss{general} grows at most linearly. Hence, it generates a globally
defined one parameter semigroup of diffeomorphisms $g_s:U_T\to
U_T$. In particular, since $g_s^*(p)=p+s$, it proves that, for all
$t_1,t_2>T$ the fibers $F_{t_1}$ and $F_{t_2}$ are diffeomorphic.

\ssec{vecf}{}
Let $p:\RR^n\to \RR$ be a polynomial map. Set 
$$
	U_t \ = \ \big\{ x\in\RR^n: \, p(x)>t\big\}.
$$
\th{vecf} There exist $T,c>0$ such that 
  \eq{vecf}
	|\n p(x)| \ > \ \frac{c}{|x|} \qquad \text{for any} \qquad x\in U_T. 
  \end{equation}
\eth
The rest of this appendix is devoted to the proof of \reft{vecf}.
\ssec{semal}{Semialgebraic sets} We will use the following results about
semialgebraic sets, cf. \cite[Appendix~A]{Hormander2}.

Recall that a subset of $\RR^n$ is called {\em semialgebraic} if it is
a finite union of finite intersection of sets defined by polynomial
equation or inequality. 

Recall also that a {\em Puiseux series} in a neighborhood of infinity
is  a series of the form
\eq{Puiseux}
	t(r) \ = \ \sum_{k=N}^\infty \, c_k r^{-\frac{k}{p}}
\end{equation}
where $p$ is an integer and $r$ is a positive variable. Here $N$ may
be a positive or negative integer or $0$.  
An important consequence of the expansion \refe{Puiseux} is that if we
choose $N$ so that $c_N\not=0$ (which is possible unless $t(r)\equiv
0$) then
$$
	t(r) \ = \ c_Nr^{-\frac{N}p}\, (1 \ + \ o(1)), \qquad
		r\ \to \infty.
$$

The proof of \reft{vecf} is based on the following theorem,
cf. \cite[Th.~A.2.6]{Hormander2}
\th{Hor}Suppose $E$ is a semialgebraic set in
   $\RR\times\RR\times\RR^n$ such that the image of the projection
   $$
	E \ \ni \ (r,t,x) \ \mapsto \ r \in \RR
   $$
   contains all large positive $r$. Then one can find Puiseux series
   $$
	t(r), \ x(r)=(x_1(r),\ldots,x_n(r))
   $$
   converging for large positive $r$ such that $(r,t(r),x(r))\in E$. 

   If 
   $$
	f(r) \ = \ \sup\big\{t: \, \text{there exists} \ x\in \RR^n \ 
		\text{such that} \ (r,t,x)\in E\big\}
   $$
   is finite and the supremum is attained for large positive $r$, one
   can take $t(r)=f(r)$.
\eth

\ssec{vecf-pr}{Proof of \reft{vecf}} Let
$E\subset\RR\times\RR\times\RR^n$ denote the set of solutions of the
following system of algebraic equations and inequalities:
\eq{syst}
	\begin{cases}
	|\n p|^2 &\le \frac1{t^2|x|^2}, \\
	p(x) &\ge t,\\
	|x|^2 &=r^2,
	\end{cases}
	\qquad 			(r,t,x)\in \RR\times\RR\times\RR^n.
\end{equation}
Set
$$
	t(r) \ = \ \sup\big\{t: \, \text{there exists} \ x\in \RR^n \ 
		\text{such that} \ (r,t,x)\in E\big\}
$$
Since the supremum is essentially taken over the compact set $\{x\in\RR^n:\,
|x|=r, x\in E\}$, it is clear that $t(r)$ is finite.

By \reft{Hor}, there exist a rational number $\alp$ and a constant $c$
such that
\eq{t(r)}
	t(r) \ = \ cr^\alp(1 \ + \ o(1)) \quad \text{as} \quad
 					t\to\infty.
\end{equation}
It also follows from \reft{Hor}, that there exist a function
$x(r)=(x_1(r)\nek x_n(r))$ defined for large $t>0$, rational numbers
$\bet_1\nek \bet_n$ and constants $c_1\nek c_n$ such that 
\eq{x(r)}
	\big(r,t(r),x(r)) \in E, \quad \text{and} \quad
	x_i(r) \ = c_i r^{\bet_i}(1+o(1)), \quad i=1\nek n.
\end{equation}
{}From \refe{syst},  we know that $|x(r)|^2=r^2$. Hence, it follows
from \refe{x(r)} that
$\bet_i\le 1$ for any $i=1\nek n$. It follows that there exists a
constant $A>0$ such that 
\eq{x'}
	\left|\frac{dx(r)}{dr}\right| \ < \ A.
\end{equation}

Suppose now that the statement of \reft{vecf} is wrong. Then $t(r)$
tends to infinity as $r\to \infty$. Hence, it follows from
\refe{t(r)}, that $\alp>0$.  Consider the function
$f(r)=p(x(r))$. Then, using the second inequality in \refe{syst}, we
obtain 
\eq{f(r)}
	\lim_{r\to\infty}f(r) \ = \  \infty.
\end{equation}
Also, from \refe{x'}, we obtain 
\eq{f'} 
	\left|\frac{df}{dr}\right| \ = \ 
	  	\left|\< \n p, \frac{dx}{dr}\>\right| \ \le \ 
			|\n p|\cdot\left| \frac{dx}{dr}\right| 
			\ < \ A|\n p|.
\end{equation}
Using \refe{f'}, \refe{t(r)} and the first inequality in \refe{syst}, we get
\eq{f'<}
	\left|\frac{df(r)}{dr}\right| \ < \ \frac{A}{r^{1+\alp}}.
\end{equation}
Since, $\alp>0$, this inequality contradicts \refe{f(r)}. $\qquad \square$


\providecommand{\bysame}{\leavevmode\hbox to3em{\hrulefill}\thinspace}

\end{document}